\documentclass{amsart}
\usepackage{amsfonts}
\usepackage{amsmath,amssymb}	
\usepackage{amsthm}
\usepackage{amscd}
\usepackage{graphics}
\usepackage{graphicx}
{
\theoremstyle{definition}
\newtheorem{Def}{Definition}
\newtheorem{Ex}{Example}
\newtheorem{Rem}{Remark}

}

\newtheorem{Thm}{Theorem}

\begin{document}
\title[Products for cobordism-like modules generic maps induce]{Products of elements of cobordism-like modules induced from generic maps}

\author{Naoki Kitazawa}
\address{Institute of Mathematics for Industry, Kyushu University, 744 Motooka, Nishi-ku Fukuoka 819-0395, Japan}
\email{n-kitazawa.imi@kyushu-u.ac.jp}
\maketitle
\begin{abstract}
Recently the author has introduced cobordism-like modules induced from generic maps whose codimensions are negative. They are generalizations of cobordism modules of manifolds. They have been
 introduced in generalizing the following theorem shown by Hiratuka and Saeki in 2013--14; for a generic map whose codimension is negative including a connected component of an inverse image of a regular value being not null-cobordant and for a space defined as all connected components of inverse images, which is a polyhedron whose dimension is equal to that of the target space, the top-dimensional homology group does not vanish. Note that such spaces are fundamental and important tools in general, in the differential topological theory of Morse functions and their higher dimensional versions and application to algebraic and differential topology of manifolds, or the global singularity theory.

In this paper, the author succeeded in defining suitable elements as the products for pairs of elements in cobordism modules which may be distinct, as in the case of the ordinary cobordism modules. They are extensions of the products for pairs of ordinary cobordism classes of manifolds.   
\end{abstract}

\section{Introduction}
\label{sec:1}
Two (oriented) closed manifolds are said to be (resp. oriented) cobordant if their disjoint union is a boundary of a compact (resp. oriented) manifold.  
The cobordism relation defined as this on the family of all closed manifolds of a fixed dimension is an equivalence relation and such relations induce a natural ring structure on all closed manifolds satisfying the following.
\begin{enumerate}
\item The sum of two elements is defined canonically by considering the disjoint union of corresponding manifolds.
\item The product of two elements is defined canonically by considering the product of corresponding manifolds.
\end{enumerate}

Such relations are fundamental and important in algebraic and differential topology of manifolds and systematic studies were started by Thom etc. (\cite{thom}). We omit more precise explanations. For precise fundamental and advanced algebraic and differential topological theory on cobordisms, see also \cite{atiyah}, \cite{stong}, \cite{milnor} etc.. We abuse terminologies on cobordisms without definitions. 

Recently in \cite{kitazawa}, the author has introduced cobordism-like modules induced from generic maps whose codimensions are negative. 
Such a module is defined by considering a weaker relation on all manifolds of a dimension equal to that of the absolute value of the codimension.; we define two (different) manifolds to be cobordant if they appear as inverse images of two regular values 
 in adjacent components of the set of all regular values.
They are generalizations of ordinary cobordism modules of manifolds.

We present the reason why the author has introduced such a module in \cite{kitazawa} with explanations on {\it Reeb spaces} and several classes of smooth maps. 

\begin{figure}
\includegraphics[width=30mm]{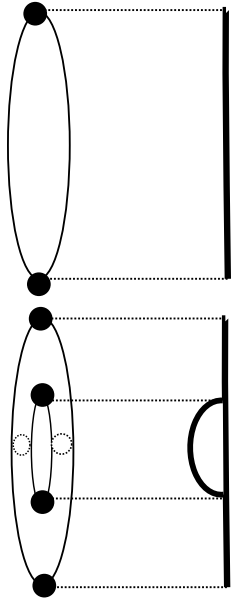}
\caption{Fundamental Morse functions and their Reeb spaces.}
\label{fig:1}
\end{figure}

First, we need to define a {\it Reeb space}. For a continuous map $c:X \rightarrow Y$, we consider the set of all connected components of inverse images of the map as a quotient space of the source space $X$ and denote it by $W_c$, call it the {\it Reeb space} of $c$ and denote the quotient map by $q_c:X \rightarrow W_c$ and a uniquely determined map by $\bar{c}$ satisfying the relation $\bar{c} \circ q_c=c$. 
FIGURE \ref{fig:1} represents explicit Reeb spaces. They represent fundamental Morse functions; a Morse function with just $2$ singular points, characterizing a $k$-dimensional homotopy sphere (except $4$-dimensional homotopy spheres not diffeomorphic to $S^4$, which are undiscovered,) topologically ($k \geq 2$), and a height function on $S^1 \times S^k$ for $k \geq 1$. Note that the Reeb space is a graph for a smooth function on a closed manifold having a finite number of singular values.

The map $c$ is said to be {\it proper} if for every compact subset $P \subset Y$, $c^{-1}(P)$ is compact.
 
Related to Reeb spaces of smooth maps whose codimensions are negative, the following has been shown.

\begin{Thm}[\cite{hiratukasaeki},\cite{hiratukasaeki2} etc.]
\label{thm:1}
For a proper smooth map $c:X \rightarrow Y$ between smooth manifolds without boundaries such that the codimension is negative and $-k$ and for suitable triangulations $k_X:K_X \rightarrow X$ and $k_Y:K_Y \rightarrow Y${\rm (}, representing canonical PL structures of $X$ and $Y${\rm )}, ${k_Y}^{-1} \circ c \circ k_X$ is a simplicial map, then the Reeb space is a polyhedron. Moreover, if there exists a connected component of an inverse image of a regular value being not null-cobordant or such that the cobordism class does not vanish, then the dimension of the Reeb space is equal to that of $Y$ and its top-dimensional homology group with coefficient ring $\mathbb{Z}/2\mathbb{Z}$ does not vanish. Furthermore, if the source and the target manifolds are oriented, then the top-dimensional homology group of the Reeb space whose coefficient ring is the $k$-dimensional smooth oriented cobordism group does not vanish.
\end{Thm}

We can apply this theorem to ({\it stable}) Morse functions, (or Morse functions such that at distinct singular points, the values are distinct,) existing densely on any smooth closed manifold, their higher dimensional versions including so-called ({\it stable}) {\it fold} maps and more general {\it stable} maps: a {\it stable} map $c$ is a smooth map such that by a slight perturbation, for the resulting map ${c}^{\prime}$, there exists a pair of diffeomorphisms $(\Phi,\phi)$ and the relation $\phi \circ c=c^{\prime} \circ \Phi$ holds or in short, a smooth map such that by a slight perturbation, the sets of all singular points and the sets of all singular values are invariant modulo diffeomorphisms. Note that stable maps exist densely in the space of all smooth maps from smooth closed manifolds into smooth manifolds without boundaries if the pairs of the dimensions of the source and the target manifolds are {\it nice}: if source and target manifolds are low-dimensional, the conditions hold. For fundamental explanations on fundamental singularity theory and geometric theory of stable maps, see \cite{golubitskyguillemin} for example.

Note also that such stable maps and Reeb spaces are fundamental and important tools in general in the differential topological theory of Morse functions and their higher dimensional versions and application to algebraic and differential topology of manifolds, or the global singularity theory. For example, Reeb spaces inherit fundamental invariants such as homology groups etc. in considerable cases. We can explicitly know this fact by FIGURE \ref{fig:1}: in the first case, the $j$-th homology group of the manifold and that of the Reeb space are isomorphic for $0 \leq j \leq k-1$ and in the second case, for $k\geq 2$, the same fact holds for $1 \leq j \leq k-1$. Moreover, Theorem \ref{thm:1} is an important result in such a branch of geometry and mathematics.

In such a stream, in \cite{kitazawa}, the author has shown an extended similar result for a presented generalized module and present related examples. That is why such a module has been introduced.

Furthermore, as an advanced study, in \cite{kitazawa3}, the author has investigated algebraic structures of such modules in cases where the codimensions of maps are $-2$.

In this paper, as another new algebraic topological study, we consider about a product of a pair of elements in these modules which may be distinct and give an answer. We define the product as an element of a module induced from a generic map whose codimension is negative and which is constructed by respecting algebraic and differential topological properties of the given maps inducing the given modules. The resulting products are regarded as extensions of products in ordinary cobordism rings. 
Throughout the present paper, let $R$ be a principle ideal domain having a unique identity element $1$ not equal to zero and for a positive integer $k$, $\mathcal{N}_k(R)$ ($\mathcal{O}_k(R)$) be a free module generated by all elements corresponding to equivalence classes obtained by considering the equivalence relation on smooth, closed and connected, (resp. oriented) manifolds defined canonically  by (resp. orientation preserving) diffeomorphisms such that distinct elements corresponding to distinct connected manifolds are mutually independent. 
Moreover, manifolds and maps between them are smooth and of class $C^{\infty}$ unless otherwise stated. The {\it singular set} of a smooth map is defined as the set of all singular points of the map, the {\it singular value set} of the map is defined as
 the image of the singular set and the {\it regular value set} of the map is the complement of the singular value set.
  
In addition, let $M$ be a closed and connected manifold of dimension $m>1$, $N$ be a connected manifold of dimension $1$ without boundary, or $\mathbb{R}$ or $S^1$, and $f$ be a map from $M$ into $N$.
   We denote the singular set of the map by $S(f)$.

The organization of the present paper is as the following. Before presenting the main theorem and the proof, we review a canonically obtained submodule $A \subset \mathcal{N}_{m-1}(R)$ (resp. $\mathcal{O}_{m-1}(R)$ oriented) compatible with a map $f:M \rightarrow N$ having a finite number of singular values based on \cite{kitazawa} and \cite{kitazawa3}. Then we show a main theorem with its proof.

The author is a member of the project and supported by the project Grant-in-Aid for Scientific Research (S) (17H06128 Principal Investigator: Osamu Saeki)
"Innovative research of geometric topology and singularities of differentiable mappings"
(https://kaken.nii.ac.jp/en/grant/KAKENHI-PROJECT-17H06128/
).
\section{Preliminaries}
\label{sec:2}

As announced in the end of the previous section, in this paper, we only treat functions whose singular value sets are finite unless otherwise stated.
For example, Morse functions, Morse-Bott functions (see for \cite{bott} for example), circle-valued versions of these functions on closed manifolds (we call {\it Morse} maps, {\it Morse-Bott} maps etc.) etc. are considered in the present paper. Reeb spaces are graphs for such functions and maps (see also \cite{masumotosaeki}, \cite{michalak}, \cite{michalak2}, \cite{sharko} etc. for such functions and maps). We introduce a fundamental object in this paper based on \cite{kitazawa} and \cite{kitazawa3}, changing terminologies to some extent.    

\begin{Def}
\label{def:1}
Let $f:M \rightarrow N$ be a map whose singular value set is finite.
For a closed interval $I$ in $N$ containing zero or one singular value such that the singular value is in the interior if it exists, let $a$ and $b$ be two distinct boundary points. For each connected component $I_c$ of the inverse image $\bar{f}^{-1}(I)$, consider all connected components of $I_c \bigcap \bar{f}^{-1}(a)$ and $I_c \bigcap \bar{f}^{-1}(b)$, which are ($m-1$)-dimensional closed and connected manifolds, and corresponding elements in $\mathcal{N}_{m-1}(R)$, and obtain the sum of all elements, respectively. We denote the sums corresponding to $a$ and $b$ by $I_{c,a}$ and $I_{c,b}$, respectively.

Let $A \subset \mathcal{N}_{m-1}(R)$ be a submodule. Let $A$ satisfy the following.   
\begin{enumerate} 
\item $A$ includes all elements represented by closed and connected manifolds not appearing as connected components of inverse images of regular values and we denote the module generated by all such elements by $A_1$.
\item For any $I_c$ before, $I_{c,a}-I_{c,b} \in A$ holds and we denote the module generated by all such elements by $A_2$.  
\item $A$ is represented as the direct sum $A_1 \oplus A_2$.
\end{enumerate}
Then $A$ is said to be {\it canonically compatible} or {\it CC} with $f$. We call $A_1$ the {\it outer part} and $A_2$ be the {\it effective part} of $A$.
If we replace $\mathcal{N}_{m-1}(R)$ by $\mathcal{O}_{m-1}(R)$, then $A$ is said to be {\it canonically oriented compatible} or {\it COC} with $f$. We define the {\it outer part} and the {\it effective part} in a same way.
\end{Def}

Note that $A$ is uniquely determined from $f$ and we denote $A$ by $CC(f)$ in the case where the manifolds are not oriented and $A$ by $COC(f)$ in the case where the manifolds are oriented. For such a module $A$, we denote $A_o$ the
 outer part and $A_e$ the effective part.
We denote the element represented by a $k$-dimensional closed (oriented) manifold $S$ by $[S] \in {\mathcal{N}}_k(R)$ (resp. $[S] \in {\mathcal{O}}_k(R)$): for a non-connected manifold, we consider the sum of the elements represented by the connected components and the coefficient of each element is positive. 

\begin{figure}
\includegraphics[width=50mm]{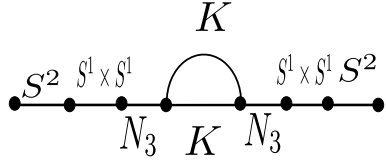}
\caption{An explicit Morse function on the product of $S^1$ and the Klein Bottle $K$.}
\label{fig:2}
\end{figure}

\begin{Ex}
For a positive integer $k>0$, let us denote the closed and non-orientable surface of genus $k$ by $N_k$ and $K:=N_1$. FIGURE \ref{fig:2} represents a fundamental Morse function on $S^1 \times K$ and manifolds represent the manifolds appearing as corresponding connected components of inverse images of regular values.  
For the module $A$ CC with the function, the effective part is generated by $[S^2]$, $[S^2]-[S^1 \times S^1]$, $[S^1 \times S^1]-[N_3]$ and $[N_3]-2[K]=[N_3]-[K \sqcup K]$.
Note that we cannot apply Theorem \ref{thm:1} immediately and can apply the generalized version in \cite{kitazawa} for the module ${\mathcal{N}}_2(R)/A$, isomorphic to $\mathbb{Z}/2\mathbb{Z}$ and regarded as the module  generated by the canonically obtained element corresponding to $[K]$ ($H_1(W_f;\mathbb{Z}/2\mathbb{Z})$ is isomorphic to $\mathbb{Z}/2\mathbb{Z}$ and not zero). 
\end{Ex}
\section{The main theorem}
\label{sec:3}
\begin{Thm}
\label{thm:2}
Let $m,m^{\prime}>1$ be integers.
Let $f:M \rightarrow N$ be a map from an $m$-dimensional closed and connected {\rm (}oriented{\rm )} manifold into a $1$-dimensional connected {\rm (}resp. oriented{\rm )} manifold with no boundary and let $f^{\prime}:M^{\prime} \rightarrow N^{\prime}$ be a
 map from an $m^{\prime}$-dimensional closed and connected {\rm (}resp. oriented{\rm )} manifold $M^{\prime}$ into a $1$-dimensional connected  {\rm (}resp. oriented{\rm )} manifold $N^{\prime}$ with no boundary.
We also assume that their singular value sets are finite.
Then we have the following.
\begin{enumerate}
\item
 There exists a map $CF(f,f^{\prime})$ on an {\rm (}$m+m^{\prime}-1${\rm )}-dimensional closed and connected {\rm (}resp. oriented{\rm )} manifold into a $1$-dimensional connected {\rm (}resp. oriented{\rm )} manifold with no boundary whose singular value set is finite and which satisfies the following.
\begin{enumerate}
\item
\label{thm:2.1}
 The module ${CC(CF(f,f^{\prime}))}_o$ {\rm (}resp. ${COC(CF(f,f^{\prime}))}_o${\rm )} is generated by all elements represented by closed and connected {\rm (}resp. oriented{\rm )} {\rm (}$m+m^{\prime}-2${\rm)}-dimensional manifolds
 we cannot represent as products of a manifold appearing as a connected component of the inverse image of a regular value of $f$ and a manifold appearing as a connected component of the inverse image of a regular value of ${f}^{\prime}$. Moreover, the module is free and distinct elements, which are all represented by closed and connected {\rm (}resp. oriented{\rm )} manifolds considered here, are independent.
\item
\label{thm:2.2}
 The module ${CC(CF(f,f^{\prime}))}_e$ {\rm (}resp. ${COC(CF(f,f^{\prime}))}_e${\rm )} is generated by the set $S$ satisfying the following.
\begin{enumerate}
\item
\label{thm:2.2.1}
 For a {\rm (}resp. an oriented{\rm )} manifold $F_M$ appearing as a connected component of the inverse image of a regular value of $f$ and for an element $[F_N]={\Sigma}_{j=1}^{k} a_k[{F_N}_k]$ of ${CC(f^{\prime})}_e$ {\rm (}resp. ${COC(f^{\prime})}_e${\rm )} represented as a linear combination of elements in $\{[{F_N}_k]\}$ where ${F_N}_k$ is a closed and connected {\rm (}resp. oriented{\rm )} manifold and where $a_k$ is an integer, ${\Sigma}_{j=1}^{k} a_k[F_M \times {F_N}_k] \in S$.
\item
\label{thm:2.2.2}
 For a {\rm (}resp. an oriented{\rm )} manifold $F_N$ appearing as a connected component of the inverse image of a regular value of $f^{\prime}$ and for an element $[F_M]={\Sigma}_{j=1}^{k} a_k[{F_M}_k]$ of ${CC(f)}_e$ {\rm (}resp. ${COC(f)}_e${\rm )} represented as a linear combination of elements in $\{[{F_M}_k]\}$ where ${F_M}_k$ is a closed and connected {\rm (}resp. oriented{\rm )} manifold and where $a_k$ is an integer, ${\Sigma}_{j=1}^{k} a_k[{F_M}_k \times {F_N}] \in S$.
\item
\label{thm:2.2.3}
 Any element of $S$ is of the form of {\rm (\ref{thm:2.2.1})} or {\rm (\ref{thm:2.2.2})}. 
\end{enumerate}
\end{enumerate}
In addition, if the $f$ and $f^{\prime}$ are Morse, then the resulting map is a Morse-Bott map.
\item Canonically, we have a homomorphism map from $$\mathcal{N}_{m-1}(R)/CC(f) \times \mathcal{N}_{{m}^{\prime}-1}(R)/CC({f}^{\prime})$$ into $$\mathcal{N}_{m+m^{\prime}-2}(R)/CF(f,f^{\prime})$$ {\rm (}resp. $$\mathcal{O}_{m-1}(R)/COC(f) \times \mathcal{O}_{{m}^{\prime}-1}(R)/COC({f}^{\prime})$$ into $$\mathcal{O}_{m+m^{\prime}-2}(R)/CF(f,f^{\prime})$${\rm )} \\
 where for two closed and connected {\rm (}resp. oriented{\rm )} manifolds, we define the resulting value as the element represented by the product of the manifolds. Moreover, the homomorphism is surjective. 
\end{enumerate}

\end{Thm}

We explain a rule about presented Reeb spaces.
If a presented Reeb space has no big dot on boundary points, then the Reeb space is for a map onto $S^1$ and the boundary points are not singular values. Manifolds represent connected components of inverse images of corresponding regular values. Note
 that there may exist empty sets in the family of the manifolds, that in such cases, Reeb spaces are empty in corresponding places and that also in such cases, we can discuss similarly.
In addition, vertical dotted lines are used to show that the Reeb spaces are for local maps, for example, for maps obtained by restrictions of the original global maps.
\begin{proof}
We show the statement in the case where manifolds are not oriented and we can show this in a same way in the case where manifolds are oriented.

The second statement on the homomorphism including the fact that the homomorphism is surjective follows from the fact that any manifold represented as a product of two manifolds, appearing as connected components of the inverse images of regular values of $f$ and $f^{\prime}$, respectively, must appear
 as a connected component of the inverse image of a regular value of the resulting map, Definition \ref{def:1} and various properties of the map in the first statement. \\ 	
\\

\noindent CASE 1 In the case where $f$ or ${f}^{\prime}$ has no singular value or where at least one of the maps gives a bundle structure over the circle. \\
\indent Let $f^{\prime}$ be such a map and let the fiber of the bundle $F$. In this case, ${CC({f}^{\prime})}_o$ is freely generated by all elements represented by closed and connected manifolds of dimension $m^{\prime}-1$ except one represented by $F$
 and ${CC({f}^{\prime})}_e$ is a trivial group by Definition \ref{def:1}. For $S$ in (\ref{thm:2.2}), it is sufficient that we cosider only (\ref{thm:2.2.2}).    
By the definitions of related terminologies, we can see that the map obtained by the composition of the canonical projection $M \times F$ onto $M$ and $f$ gives the desired map. We can see also that the other facts we should show is true.

See also FIGURE \ref{fig:1} for explicit construction.

\begin{figure}
\includegraphics[width=50mm]{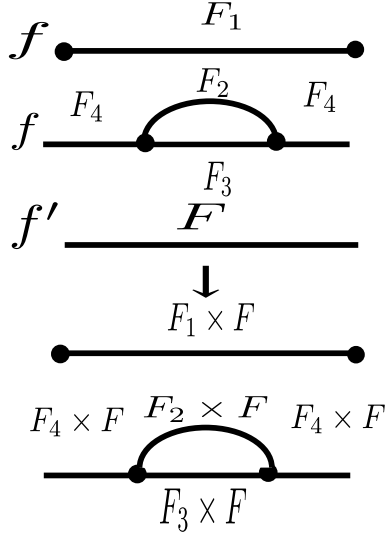}
\caption{Construction of desired maps for CASE 1.}
\label{fig:3}
\end{figure}

\ \\

\noindent CASE 2  In the case where $f$ and ${f}^{\prime}$ have singular values. \\
\indent For each singular value of $f$ and $f^{\prime}$, we consider a closed small interval such that the singular value is in the interior and consider connected components of the inverse image of the interval.
For each singular value of $f$, the closed small interval and each connected component of the inverse image of the interval, we change the map into one as in the upper part of FIGURE \ref{fig:4}. We change the local map into a map obtained
 by three copies of the original local map canonically. For each singular value of ${f}^{\prime}$, the closed small interval and each connected component of the inverse image of the interval, we consider the map as in the lower part of FIGURE \ref{fig:4}. The map is obtained by gluing two copies of the original local map canonically. 
\begin{figure}
\includegraphics[width=50mm]{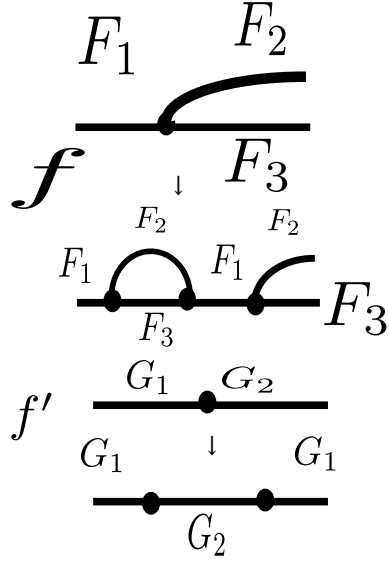}

\caption{Local changes of the maps around each connected component of the inverse image of a small closed interval including a singular value in the interior.}
\label{fig:4}
\end{figure}
Next, for $f$, we change the obtained local map again through FIGURE \ref{fig:5} and FIGURE \ref{fig:6} for example and we can naturally generalize the discussion here and we will explain again. 
\begin{figure}
\includegraphics[width=50mm]{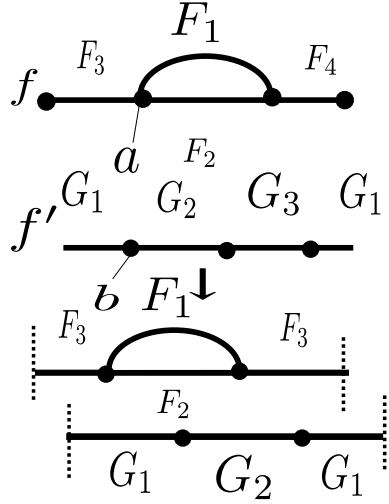}
\caption{Local maps around vertices $a$ and $b$, corresponding to singular values.}
\label{fig:5}
\end{figure}
\begin{figure}
\includegraphics[width=50mm]{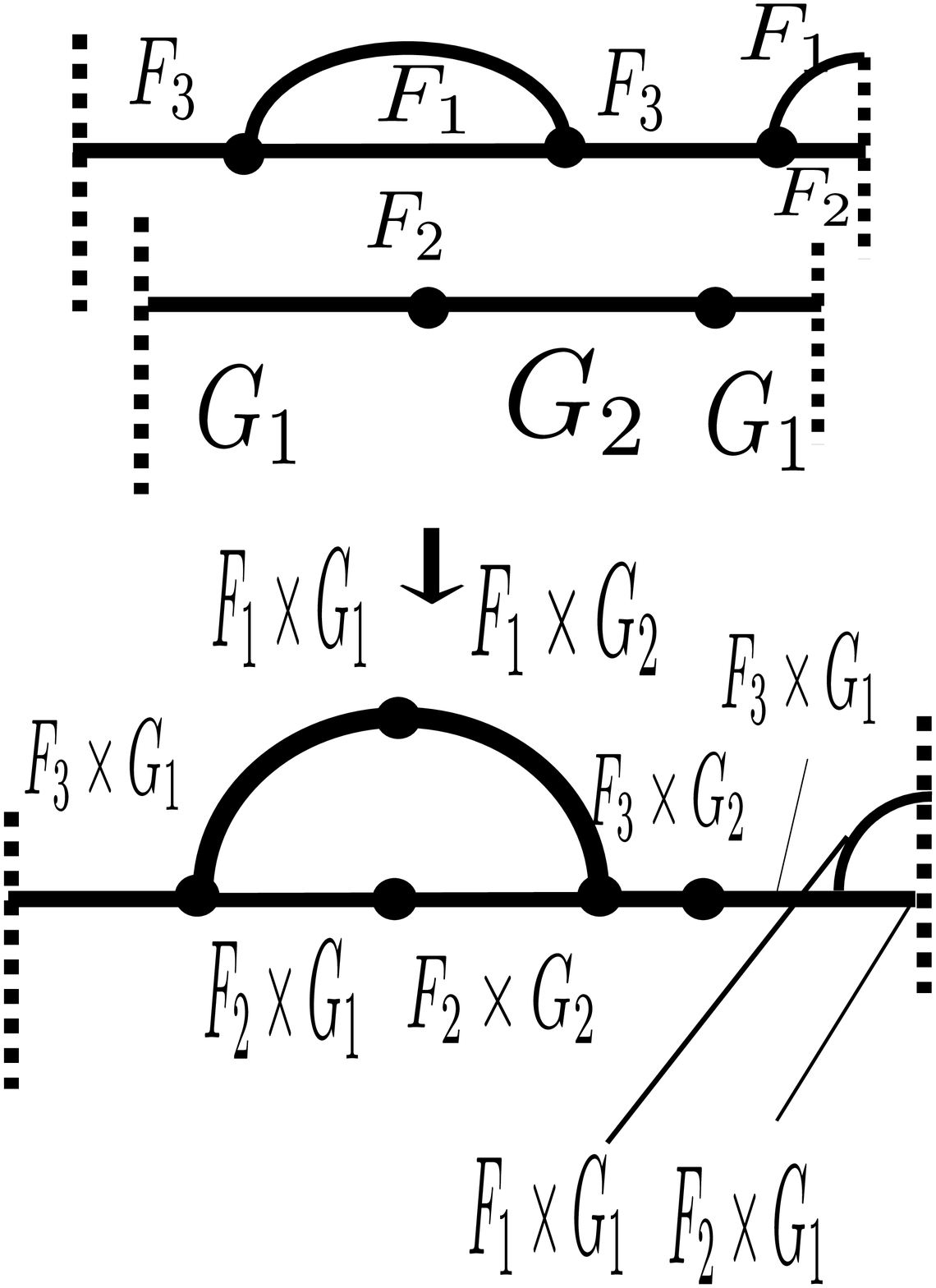}

\caption{Local changes of the maps in FIGURE \ref{fig:5}.}
\label{fig:6}
\end{figure}

For a singular value $a$ of $f$ and a singular value $b$ of $f^{\prime}$, we consider a closed small interval $[a_1,a_2] \subset \mathbb{R} \subset N$ such that the singular value $a$ is in the interior and the intersection $I_{c,i}$ of $f^{-1}(a_i)$ and a connected component of the inverse image of the interval for $i=1,2$ and consider also a closed small interval $[b_1,b_2] \subset \mathbb{R} \subset N^{\prime}$ such that the singular value $b$ is in the interior and the intersection $I_{c^{\prime},i}$ of ${f^{\prime}}^{-1}(b_i)$ and a connected component of the inverse image of the interval for $i=1,2$. We can change the local map into a local map such that the inverse image of a regular value is diffeomorphic to $I_{c,1} \times I_{c^{\prime},1}$, $I_{c,2} \times I_{c^{\prime},1}$,
 $I_{c,2} \times I_{c^{\prime},2}$, $I_{c,1} \times I_{c^{\prime},2}$, $I_{c,1} \times I_{c^{\prime},1}$ or $I_{c,2} \times I_{c^{\prime},1}$ and the inverse image of a regular value changes in order as the value increases from $a_1$ to $a_2$. Note that locally around a singular value, the form of the map is as presented in CASE 1 or represented as the composition of a canonical projection and the original map. 

For a fixed singular value of $f^{\prime}$, a closed small interval including the value in the interior and a fixed connected component of the inverse image of the interval, by doing such a procedure for each singular value of $f$, a closed small interval including the value in the interior and each connected component of the inverse image of the interval, we obtain a map into $N$ whose singular value set is finite.

For each singular value of $f^{\prime}$, a closed small interval including the value in the interior and each connected component of the inverse image of the interval, we can do such a procedure and as a result, we obtain a map. \\
\indent In the previous step, we obtain a map whose singular value set is finite. If the resulting source manifold is not connected, then by deforming the map by a suitable deformation by a family of isotopies and using a technique also used in \cite{kitazawa3} in cases where maps are of codimension $-2$, in \cite{michalak} in cases where inverse images of regular values are disjoint unions of standard spheres for example etc., we can make the source manifold connected preserving the module CC with the function. In other words, we change a local map on the inverse image of a closed interval in the regular value set into a local Morse function with just $1$ singular value. We can do this by virtue of a fundamental fact on handle attachments corresponding to singular points. See also FIGURE \ref{fig:7} and FIGURE \ref{fig:8}.  

\begin{figure}
\includegraphics[width=50mm]{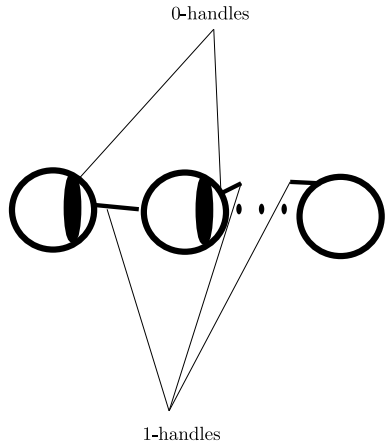}

\caption{Attaching $0$-handles and $1$-handles at once.}
\label{fig:7}
\end{figure}

\begin{figure}
\includegraphics[width=50mm]{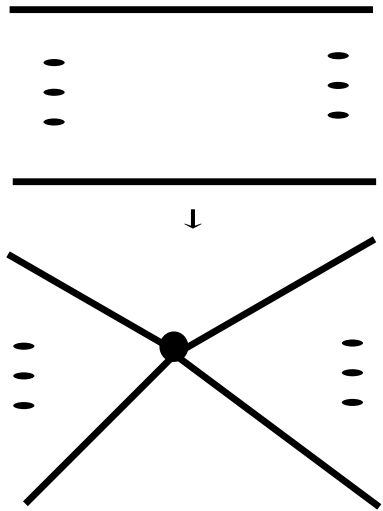}

\caption{The Reeb space of a local Morse function corresponding to the attachment of handles in FIGURE \ref{fig:7}.}
\label{fig:8}
\end{figure}

This gives a desired map $CF(f,f^{\prime})$.  \\
\\
\indent Last, if the original maps are Morse, then by the method of construction, the resulting map is a Morse-Bott map. This completes the proof.
\end{proof}

\begin{Ex}
\label{ex:2}
\begin{enumerate}
\item
Let $f$ be a Morse function on an $m$-dimensional closed and connected manifold whose Reeb space is as FIGURE \ref{fig:9} where $\Sigma$ is a homotopy sphere of dimension $m-1$. For such a Morse function, see also \cite{kitazawa2}. 
Let $f^{\prime}$ give a bundle over the circle whose fiber is an ($m^{\prime}-1$)-dimensional closed and connected manifold $F$. We demonstrate the construction of the proof of CASE 1 of Theorem \ref{thm:2}.
\begin{figure}
\includegraphics[width=50mm]{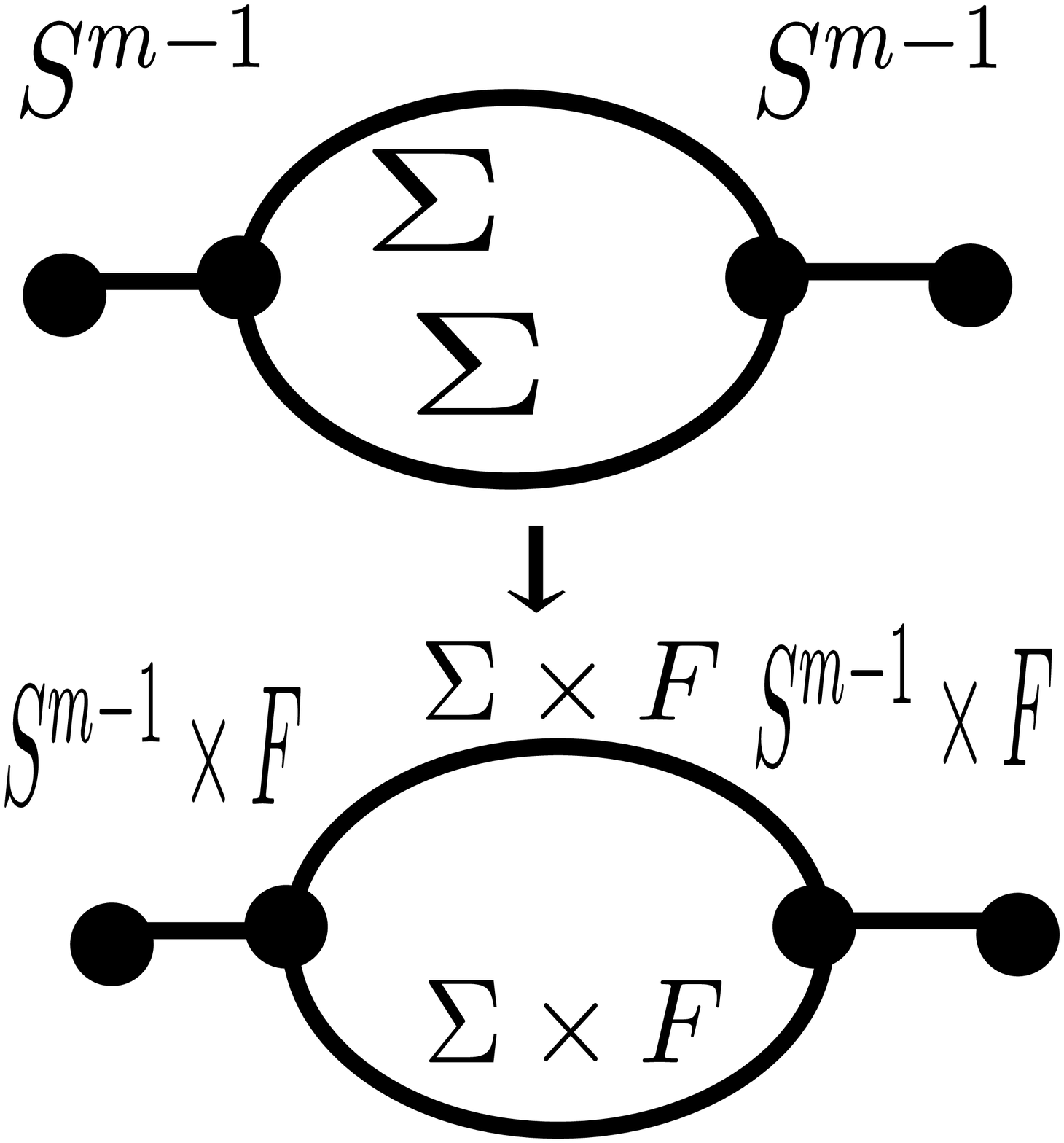}

\caption{An example to which we can apply CASE 1 of the proof of Theorem \ref{thm:2}.}
\label{fig:9}
\end{figure}
\item
Let $f$ be a Morse function on an $m$-dimensional closed and connected manifold whose Reeb space is as FIGURE \ref{fig:10}, which is same as the map of the previous example. Let $f^{\prime}$ be a Morse function with just $2$-singular points, characterizing $m^{\prime}$-dimensional homotopy sphere topologically. Through FIGUREs \ref{fig:10}--\ref{fig:13}, we demonstrate the construction of CASE 2 of the proof of Theorem \ref{thm:2}. Note that dotted lines
 along the Reeb spaces mean that the inverse images are empty (dotted lines are in fact not in the Reeb spaces).
\begin{figure}
\includegraphics[width=50mm]{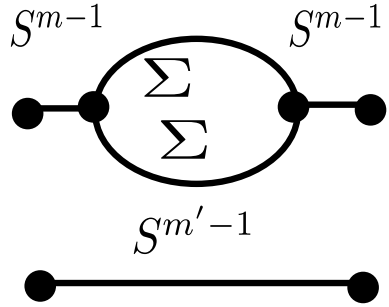}

\caption{An example to which we can apply CASE 2 of the proof of Theorem \ref{thm:2}.}
\label{fig:10}
\end{figure}
\begin{figure}
\includegraphics[width=50mm]{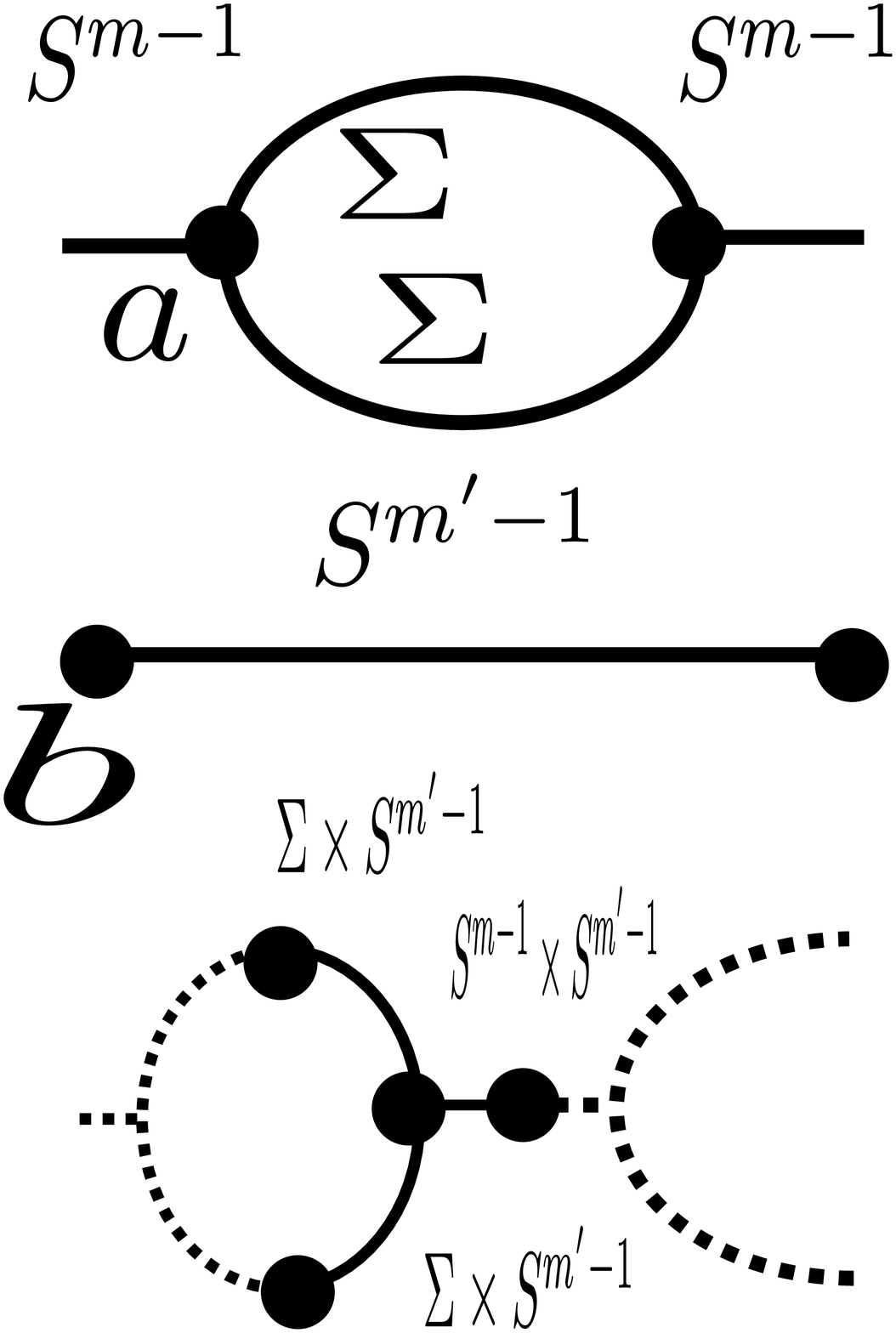}

\caption{In the case of FIGURE \ref{fig:10}, for vertices $a$ and $b$, we apply the construction.}
\label{fig:11}
\end{figure}

\begin{figure}
\includegraphics[width=50mm]{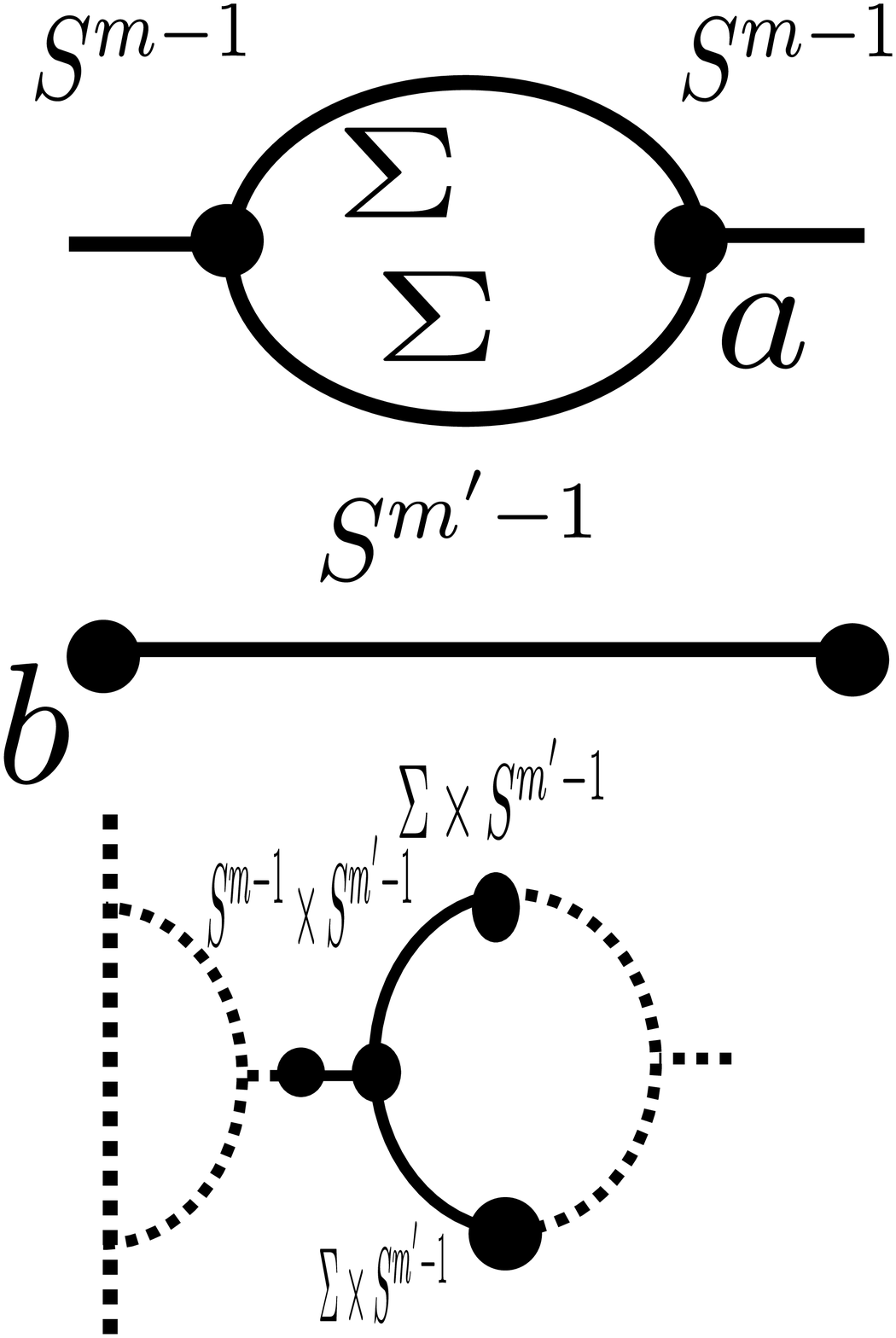}

\caption{In the case of FIGURE \ref{fig:10}, for vertices $a$ and $b$, we apply the construction.}
\label{fig:12}
\end{figure}
\begin{figure}
\includegraphics[width=50mm]{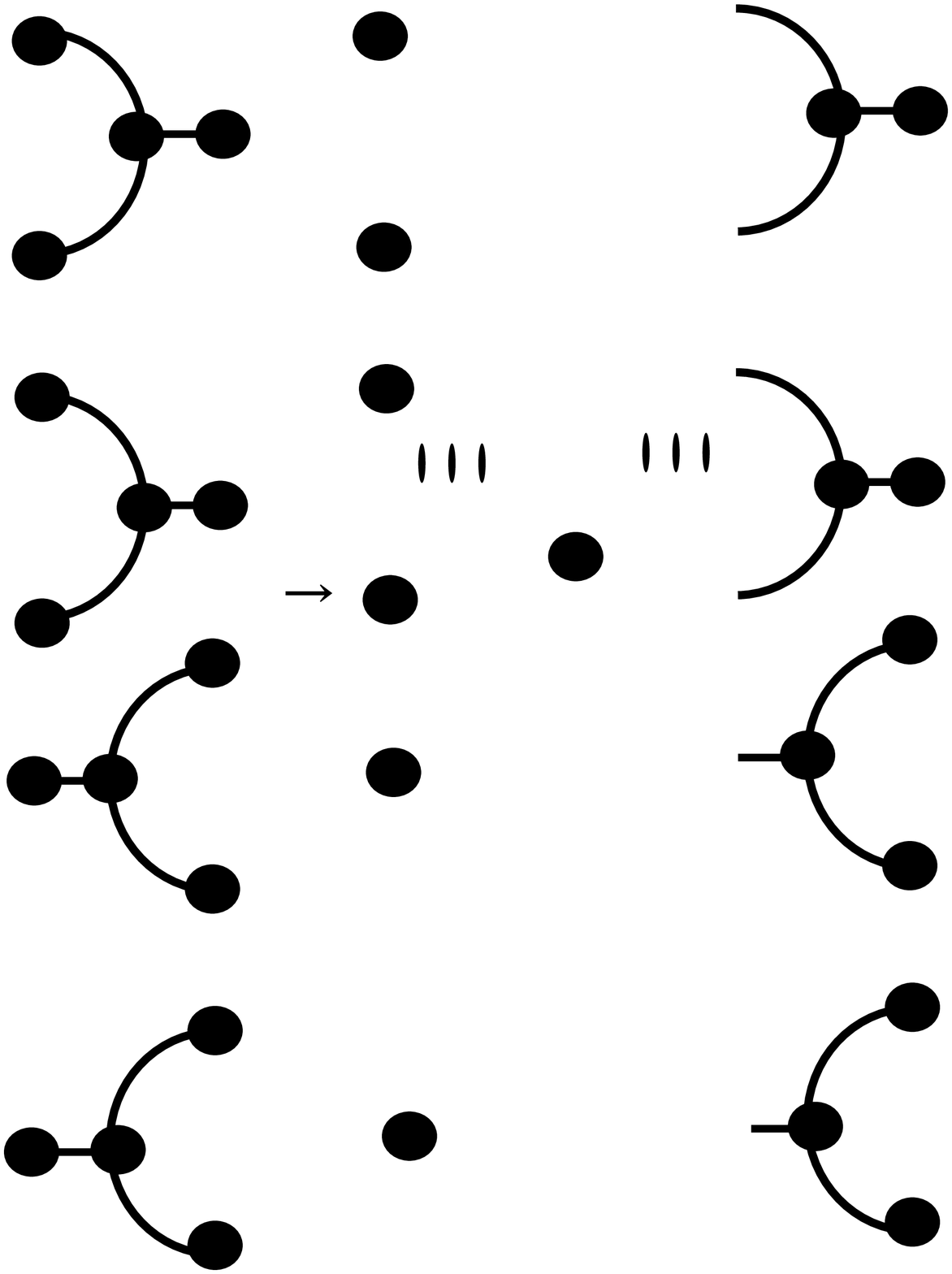}

\caption{The last step to make the source manifold connected in the case of FIGURE \ref{fig:10} (we omit information of manifolds appearing as connected components of inverse images of regular values).}
\label{fig:13}
\end{figure}
\end{enumerate}
\end{Ex}
\begin{Rem}
In both cases of Example \ref{ex:2}, the manifolds appearing as connected components of inverse images of regular values can be unique modulo diffeomorphisms for suitable pairs $(m,m^{\prime})$ of dimensions. Even if $\Sigma$ is exotic (in this case the module $\mathcal{N}_{m-1}(\mathbb{Z})/CC(f)$ is isomorphic to $\mathbb{Z}/2\mathbb{Z}$), then the
 manifolds are always diffeomorphic to $S^{m-1} \times S^{m^{\prime}-1}$ in suitable cases: we can make such a case by choosing $m$ and $m^{\prime}$ so that the number of connected components of the diffeomorphism group of $S^{m+m^{\prime}-2}$ is one: for example, in the case where the relation $m+m^{\prime}-2=11$ holds (see \cite{crowleyschick} and \cite{kitazawa3} for example). 

In the first case, if the number of the connected components of the diffeomorphism group of $S^{m+m^{\prime}-2}$ is small, then the homomorphism of the second statement of Theorem \ref{thm:2} seems to be not injective (and the zero map) in considerable cases.
\end{Rem}
\begin{Rem}
In considerable cases, it seems to be true that we can construct desired maps in Theorem \ref{thm:2} as ones simpler than ones obtained by the construction of the proof. For example, in the second example of Example \ref{ex:2}, maps obtained in FIGURE \ref{fig:11} and FIGURE \ref{fig:12} are also answers.  
\end{Rem}


\begin{thebibliography}{25} 
\bibitem{atiyah} M. F. Atiyah, \textsl{Bordism and cobordism}, Proc. Camb. Phil. Sco. 57 (1961)., 200--208.
\bibitem{bott} R. Bott, \textsl{Nondegenerate critical manifolds}, Ann. of Math. 60 (1954), 248--261.
\bibitem{crowleyschick} D. Crowley, T. Schick, \textsl{The Gromoll filtration, KO-chracteristic classes and metrics of positive scalar curvature}, Geometry  Topoogy 17 (3) (2013), 1773--1789, arxiv:1204.6474.
\bibitem{golubitskyguillemin} M. Golubitsky and V. Guillemin, \textsl{Stable Mappings and Their Singularities}, Graduate Texts in Mathematics (14), Springer-Verlag(1974).
\bibitem{hiratukasaeki} J. T. Hiratuka and O. Saeki, \textsl{Triangulating Stein factorizations of generic maps and Euler Characteristic formulas}, RIMS Kokyuroku Bessatsu B38 (2013), 61--89. 
\bibitem{hiratukasaeki2} J. T. Hiratuka and O. Saeki, \textsl{Connected components of regular fibers of differentiable maps}, in "Topics on Real and Complex Singularities", Proceedings of the 4th Japanese-Australian Workshop (JARCS4), Kobe 2011,  World Scientific, 2014, 61--73. 
\bibitem{kitazawa} N. Kitazawa, \textsl{Smooth maps compatible with simplicialstructures and inverse images}, submitted to a refereed journal, arxiv:1802.06381 (arxiv:1802.06381v6).
\bibitem{kitazawa2} N. Kitazawa, \textsl{Lifts of spherical Morse functions}, submitted to a refereed journal, arxiv:1805.05852 (arxiv:1805.05852v3).
\bibitem{kitazawa3} N. Kitazawa, \textsl{Structures of cobordism-like modules induced from generic maps of codimension -2}, submitted to a refereed journal, arxiv:1901.04994 (arxiv:1901.04994v1).
\bibitem{masumotosaeki} Y. Masumoto and O. Saeki, \textsl{A smooth function on a manifold with given Reeb graph}, Kyushu J. Math. 65 (2011), 75--84.
\bibitem{michalak} L. P. Michalak, \textsl{Realization of a graph as the Reeb graph of a Morse function on a manifold}. to appear in Topol. Methods Nonlinear Anal., Advance publication (2018), 14pp, arxiv:1805.06727.
\bibitem{michalak2} L. P. Michalak, \textsl{Combinatorial modifications of Reeb graphs and the realization problem}, arxiv:1811.08031.
\bibitem{milnor} J. Milnor, \textsl{Lectures on the h-cobordism theorem}, Math. Notes, Princeton Univ. Press, Princeton, N.J. 1965.
\bibitem{reeb} G. Reeb, \textsl{Sur les points singuliers d\'{}une forme de Pfaff compl\'{e}tement int\`{e}grable ou d\'{}une fonction num\'{e}rique}, Comptes Rendus
 Hebdomadaires des S\'{e}ances de I\'{}Acad\'{e}mie des Sciences 222 (1946), 847--849.
\bibitem{sharko} V. Sharko, \textsl{About Kronrod-Reeb graph of a function on a manifold}, Methods of Functional Analysis and
\bibitem{shiota} M. Shiota, \textsl{Thom's conjecture on triangulations of maps}, Topology 39 (2000), 383--399. 
\bibitem{stong} R. E. Stong, \textsl{Notes on cobordism theory}, Princeton Universty Press, 1968.
\bibitem{thom} R. Thom, \textsl{Quelques propri\'{e}t\'{e}s globales des vari\'{e}t\'{e}s diff\'{e}rentiables}, Commentarii Mathematici Helvetici 28, 17--86 (1954).

\end{thebibliography}
\end{document}